\newcommand{\al}{\alpha}               
\newcommand{\ga}{\gamma}               \newcommand{\Ga}{\Gamma}

\newcommand{\veps}{\varepsilon}        \newcommand{\vphi}{\varphi}

\newcommand{\cal}{\mathcal}
           
           \newcommand{\calf}{{\cal F}}
           
\newcommand{\calm}{{\cal M}}

       \newcommand{\Dom}{{\rm Dom}}
\newcommand{\Fix}{{\rm Fix}}           \newcommand{\Pos}{{\rm Pos}}
    
\newcommand{\incl}{\subseteq}          
\newcommand{\es}{\emptyset}            \newcommand{\sm}{\setminus}
\newcommand{\impl}{\Rightarrow}        \newcommand{\limpl}{\Longrightarrow}
    
\newcommand{\oo}{\infty}

                 \newcommand{\sk}{\smallskip}
           \newcommand{\n}{\noindent}

                \def\R+oo{R_+\cup\{\oo\}}

\newcommand{\barr}{\begin{array}}          \newcommand{\earr}{\end{array}}
\newcommand{\beq}{\begin{equation}}        \newcommand{\eeq}{\end{equation}}
\newcommand{\bit}{\begin{itemize}}         \newcommand{\eit}{\end{itemize}}
\newcommand{\blemma}{\begin{lemma}}        \newcommand{\elemma}{\end{lemma}}
\newcommand{\bprop}{\begin{proposition}}   \newcommand{\eprop}{\end{proposition}}
\newcommand{\bproof}{\begin{proof}}        \newcommand{\eproof}{\end{proof}}
\newcommand{\brem}{\begin{remark}}   
\newcommand{\erem}{\end{remark}}
\newcommand{\btab}{\begin{tabular}}        \newcommand{\etab}{\end{tabular}}
\newcommand{\btheorem}{\begin{theorem}}    \newcommand{\etheorem}{\end{theorem}}

\documentclass[reqno]{amsart}

\newtheorem{theorem}{\bf Theorem}

\newtheorem{lemma}{\bf Lemma}
\newtheorem{proposition}{\bf Proposition}
\newtheorem{remark}{\bf Remark}

\begin{document}

\title
[Implicit Contractive Maps in Ordered Metric Spaces]
{IMPLICIT CONTRACTIVE MAPS \\
IN ORDERED METRIC SPACES}

\author{Mihai Turinici}
\address{
"A. Myller" Mathematical Seminar;
"A. I. Cuza" University;
700506 Ia\c{s}i, Romania
}
\email{mturi@uaic.ro}


\subjclass[2010]{
47H10 (Primary), 54H25 (Secondary).
}

\keywords{
Ordered metric space, increasing operator,
fixed point, compatible normal and positive function, 
implicit and explicit contraction.
}

\begin{abstract}
Further extensions are given to the 
fixed point result (for implicit contractions)
due to Altun and Simsek
[Fixed Point Th. Appl., Volume 2010,
Article ID 621469].
Some connections with 
related statements in the area due to
Agarwal, El-Gebeily and O'Regan
[Appl. Anal., 87 (2008), 109-116]
are also discussed.
Finally, the old approach in Turinici 
[An. \c{S}t. Univ. "A. I. Cuza" Ia\c{s}i, 22 (1976), 177-180]
is presented, for historical reasons.
\end{abstract}

\maketitle

\section{Introduction}
\setcounter{equation}{0}

Let $X$ be a nonempty set,
$d:X\times X\to R_+:=[0,\oo[$ be a metric on $X$ and
$(\le)$ be a {\it quasi-order} 
(i.e.: reflexive transitive relation) over it;
the resulting triple 
$(X,d,\le)$ will be referred to as a 
{\it quasi-ordered metric space}.
Further, take some $T\in \calf(X)$.
[Here, given the nonempty sets $A$  and $B$, 
$\calf(A,B)$ stands for 
the class of all functions $f:A\to B$; 
when $A=B$, we write $\calf(A,A)$ as $\calf(A)$].
The basic conditions to be posed upon these data are
\bit
\item[(a01)]
$(X,d,\le)$ is complete 
(each ascending $d$-Cauchy sequence is $d$-convergent)
\item[(a02)]
$X(T,\le):=\{x\in X; x\le Tx\}$ is nonempty
\item[(a03)]
$T$ is $(\le)$-increasing ($x\le y$ implies $Tx\le Ty$).
\eit
Denote $\Fix(T):=\{z\in X; z=Tz\}$; any point of it will be called
{\it fixed} under $T$.
These are to be determined in the context below
(cf. Rus \cite[Ch 2, Sect 2.2]{rus-2001}):

{\bf 1a)}
We say that $x\in X(T,\le)$ is a {\it Picard point} (modulo $(d,\le;T)$) when 
$(T^nx)$ converges and $\lim_n T^nx$ is in $\Fix(T)$

{\bf 1b)}
If this holds for each $x\in X(T,\le)$,
we say that $T$ is a {\it Picard operator} (modulo $(d,\le)$);
and, if in addition, 
$\Fix(T)$ is $(\le)$-{\it singleton} [$z,w\in \Fix(T)$ and $z\le w$ limply $z=w$],
then $T$ is called a {\it global Picard operator} (modulo $(d,\le)$). 
\sk

Let $\calf(in)(R_+)$ stand for  
the subclass of all {\it increasing} $\vphi\in \calf(R_+)$;
and $\calf(re)(R_+)$ be the subclass of all
$\vphi\in \calf(R_+)$ with the (strong) {\it regressive} property:
[$\vphi(0)=0$; $\vphi(t)< t$, $\forall t> 0$].
We say that 
$\vphi\in \calf(in,re)(R_+):=\vphi\in \calf(in)(R_+)\cap \calf(re)(R_+)$ 
is a {\it Matkowski function} if
[$\vphi^n(t)\to 0$, for all $t> 0$].

Given $\vphi\in \calf(R_+)$, call the self-map $T$, 
{\it $(d,\le;\vphi)$-contractive}, if
\bit
\item[(a04)]
$d(Tx,Ty)\le \vphi(d(x,y))$, $\forall x,y\in X$, $x\le y$.
\eit
The following answer to the posed question is available.

\btheorem \label{t1}
Suppose that (in addition to (a01)-(a03))
one of the conditions below is fulfilled
\bit
\item[(a05)]
$T$ is $(\le)$-continuous: \\
$(x_n)$ is ascending and $x_n\to x$ imply $Tx_n\to Tx$
\item[(a06)]
$(\le)$ is self-closed: \\
$(x_n)$=ascending and $x_n\to x$ imply $x_n\le x$, $\forall n$.
\eit
Further, let the selfmap $T$ be $(d,\le;\vphi)$-contractive,
where $\vphi\in \calf(in,re)(R_+)$ is a Matkowski function. 
Then, $T$ is a global Picard operator (modulo $(d,\le)$).
\etheorem

This result was obtained in 1986 by
Turinici \cite{turinici-1986}
(over the class of ordered metrizable uniform spaces).
Note that, in the amorphous case
$(\le)=X\times X$, 
(a06) is fulfilled and (a01) becomes
\bit
\item[(a07)]
$(X,d)$ is complete 
(each $d$-Cauchy sequence is $d$-convergent);
\eit
the corresponding version of Theorem \ref{t1} 
is nothing else than the  1975 statement in
Matkowski \cite{matkowski-1975},
comparable with the one in
Boyd and Wong \cite{boyd-wong-1969}.
On the other hand, when $(\le)$ is an order and (a07) 
is again holding, Theorem \ref{t1} is just the 2008
statement due to
Agarwal et al \cite[Theorem 2.1]{agarwal-el-gebeily-o-regan-2008}; cf. 
O'Regan and Petru\c{s}el \cite{o-regan-petrusel-2008}.
In particular, when $\vphi$ is {\it linear}
($\vphi(t)=\al t$, $t\in R_+$, for some $\al\in [0,1[$)
this version of Theorem \ref{t1} gives the 
statement in 
Ran and Reurings \cite{ran-reurings-2004};
see also 
Nieto and Rodriguez-Lopez \cite{nieto-rodriguez-lopez-2005}. 

The obtained variants of Theorem \ref{t1} 
found some useful applications to existence theorems 
for linear and nonlinear operator equations; 
see the quoted papers for details.
As a consequence, the question of extending this result in 
useful from both a theoretical and a practical perspective.
Some interesting results of this type were proposed,
in the amorphous case, by
Wardowski \cite{wardowski-2012},
via contractive conditions 
\bit
\item[(a08)]
$F(d(Tx,Ty),d(x,y))\le 0$, $\forall x,y\in X$;
\eit
where $F:R_+^2\to R$ is an appropriate function.
However, as shown in  
Turinici \cite{turinici-2012},
all these are reducible to 
Matkowski's 
\cite{matkowski-1975}.
An extended version of them 
was given (in the ordered framework) by
Altun and Simsek \cite{altun-simsek-2010},
by means of contractive conditions like 
\bit
\item[(a09)]
$F(d(Tx,Ty),d(x,y),d(x,Tx),d(y,Ty),d(x,Ty),d(Tx,y))\le 0$,\\
for all $x,y\in X$ with $x\le y$;
\eit
where $F:R_+^6\to R$ is a function.
However, it does not include
in a complete manner the "explicit" result 
(comparable with Theorem \ref{t1})
due to 
Agarwal et al 
\cite[Theorem 2.2]{agarwal-el-gebeily-o-regan-2008};
so, we may ask whether this is removable.
It is our aim in the present exposition to state (in Section 3) 
a further extension of this implicit fixed point principle 
which includes in a complete manner 
(cf. Section 5)
the explicit 2008 result above,
as well as (according to Section 4), Theorem \ref{t1} itself.
The preliminary facts for these developments are given in 
Section 2.

Finally, note that in almost all papers 
based on implicit techniques --
including the ones in 
Akkouchi \cite{akkouchi-2011}
or 
Berinde and Vetro \cite{berinde-vetro-2012}
(see also 
Nashine et al \cite{nashine-kadelburg-kumam-2012}) --
it is asserted that the starting point 
in the area is represented by the 
contributions due to 
Popa 
\cite{popa-1997}, 
\cite{popa-1999},
\cite{popa-2003}.
Unfortunately, all these affirmations are false; 
to verify our claim, we present, in Section 6,
an "old" implicit approach 
-- obtained four decades ago --
by Turinici \cite{turinici-1976}.
Further aspects will be discussed elsewhere.

\section{Preliminaries}
\setcounter{equation}{0}

In the following, some preliminary facts involving 
real functions and (standard) metric spaces are given.
\sk

{\bf (A)}
Let $F\in\calf(R_+^6,R)$ be a function. 

{\bf 2a)}
Call it {\it compatible}, provided:
\bit
\item[(b01)]
for each couple of sequences 
$(r_n; n\ge 0)$ in $R_+^0:=]0,\oo[$  \\ 
and  $(s_n; n\ge 0)$ in $R_+$ with 
$F(r_n,r_{n-1},r_{n-1},r_n,s_{n-1},0)\le 0$, $\forall n\ge 1$
and  $|s_{n-1}-r_{n-1}|\le r_{n}$, $\forall n\ge 1$, 
we must have $r_n\to 0$ (hence, $s_n\to 0$).
\eit

{\bf 2b)}
Further, let us say that $F\in \calf(R_+^6,R)$
is {\it (3,4)-normal}, in case
\bit
\item[(b02)]
$F(r,r,0,0,r,r)> 0$,\ for all $r>0$.
\eit

The next property will necessitate some conventions.
Take some point (in $R_+^6$) $W=(w_1,w_2,w_3,w_4,w_5,w_6)$;
as well as a rank $j\in \{1,2,3,4,5,6\}$.
We say that the sequence 
$(t^n:=(t_1^n,t_2^n,t_3^n,t_4^n,t_5^n,t_6^n); n\ge 0)$ in $R_+^6$ is

{\bf 2c)}
{\it j-right} at $W$, if
$t_i^n\to w_i$, $i\ne j$, $t_j^n\to w_j+$.

{\bf 2d)}
{\it j-point} at $W$, if
$t_i^n\to w_i$, $i\ne j$, $t_j^n=w_j$, $\forall n$.

\n
[Here, $z_n\to z+$ means: $z_n\to z$ and
$z_n> z$, $\forall n$].

{\bf 2e)}
Given $b> 0$, call the function $F$,
{\it 2-right-lim-positive} at $b$,
if 
$\limsup_n F(t^n)> 0$, 
for each 2-right at $(b,b,0,0,b,b)$, 
sequence $(t^n; n\ge 0)$ in $(R_+^0)^6$.
The class of all these $b> 0$ will be denoted as
$\Pos(2-right-lim;F)$.
In this case, we say that 
$F$ is {\it almost 2-right-lim-positive}, if
$\Theta:=\Pos(2-right-lim;F)$ is $(>)$-cofinal in 
$R_+^0$ [for each $\veps\in R_+^0$ there exists 
$\theta\in \Theta$ with $\veps> \theta$].

{\bf 2f)}
Given $b> 0$, call the function $F$,
{\it 4-point-lim-positive} at $b$, 
if 
$\limsup_n F(t^n)> 0$, 
for each 4-point at $(b,0,0,b,b,0)$ sequence 
$(t^n; n\ge 0)$ in $(R_+^0)^6$.
The class of all these $b> 0$ will be denoted as
$\Pos(4-point-lim;F)$.
In this case, we say that 
$F$ is {\it 4-point-lim-positive}, if
this last set is 
identical with $R_+^0$.
\sk

{\bf (B)}
We are now passing to another fact.
Let $(X,d)$ be a metric space.
Call the sequence $(x_n; n\ge 0)$ in $X$, 
{\it $d$-semi-Cauchy}, provided 
$d(x_n,x_{n+1})\to 0$.
In the following, a useful property is described for such
sequences which are not $d$-Cauchy.

\bprop \label{p1}
Suppose that $(x_n; n\ge 0)$ is a sequence in $X$ with
\bit
\item[(b03)] 
$r_n:=d(x_n,x_{n+1})> 0$, for all $n$
\item[(b04)]
$(x_n; n\ge 0)$ is $d$-semi-Cauchy but not $d$-Cauchy.
\eit
Further, let $\Theta$ be be a $(>)$-cofinal part of $R_+^0$.
There exist then a number $b\in \Theta$, a rank $j(b)\ge 0$, 
and a couple
of rank-sequences $(m(j); j\ge 0)$, $(n(j); j\ge 0)$, with
\beq \label{201}
j\le m(j)< n(j),\ d(x_{m(j)},x_{n(j)})> b,\ 
\forall j\ge 0
\eeq
\beq \label{202}
n(j)-m(j)\ge 2,\ d(x_{m(j)},x_{n(j)-1})\le b,\ 
\forall j\ge j(b)
\eeq
\beq \label{203}
\mbox{  
\btab{l}
$(u_j(0,0):=d(x_{m(j)},x_{n(j)}); j\ge 0)$ 
is a sequence in $R_+^0$ \\
with\ $u_j(0,0)\to b+$ as $j\to \oo$
\etab
}
\eeq
\beq \label{204}
\mbox{ 
\btab{l}
$(u_j(p,q):=d(x_{m(j)+p},x_{n(j)+q}); j\ge j(b))$ 
is a sequence in $R_+^0$ \\
with $u_j(p,q)\to b$ as $j\to \oo$,\ $\forall p,q\in \{0,1\}$. 
\etab
}
\eeq
\eprop

\bproof
By definition, the $d$-Cauchy property of our sequence writes:
\bit
\item[]
$\forall \veps\in R_+^0$, $\exists k=k(\veps)$:\
$k\le m< n$ $\limpl$ $d(x_m,x_n)\le \veps$.
\eit
As $\Theta\incl R_+^0$ is $(>)$-cofinal in $R_+^0$, 
this property may be also written as
\bit
\item[]
$\forall \theta\in \Theta$, $\exists k=k(\theta)$:\
$k\le m< n$ $\limpl$ $d(x_m,x_n)\le \theta$.
\eit
The negation of this property means:
there exists $b\in \Theta$
such that, $\forall j\ge 0$:
\beq \label{205}
A(j):=\{(m,n)\in N\times N; j\le m< n, d(x_m,x_n)> b\}
\ne \es.
\eeq
Having this precise, denote, for each $j\ge 0$,
\bit
\item[(b05)]
$m(j)=\min \Dom(A(j))$,\ $n(j)=\min A(m(j))$.
\eit
The couple of rank-sequences 
$(m(j); j\ge 0)$, $(n(j); j\ge 0)$  fulfills (\ref{201});
hence, the first half of (\ref{203}). 
On the other hand, letting
$j(b)$ be such that
\beq \label{206}
\mbox{
$r_i:=d(x_i,x_{i+1})< b/3$,\ for all $i\ge j(b)$,
}
\eeq
it is clear that  (\ref{202}) holds too.
This in turn yields, $\forall j\ge j(b)$;
$$  
b< d(x_{m(j)},x_{n(j)})\le 
d(x_{m(j)},x_{n(j)-1})+r_{n(j)-1}\le  \\
b+r_{n(j)-1};
$$
so, passing to limit as $j\to \oo$ gives 
the second half of (\ref{203}).
Finally, $\forall j\ge j(b)$,
$$
d(x_{m(j)},x_{n(j)+1})\le 
d(x_{m(j)},x_{n(j)})+r_{n(j)},
$$
$$
d(x_{m(j)},x_{n(j)+1})\ge 
d(x_{m(j)},x_{n(j)})-r_{n(j)}> 2b/3.
$$
This gives the case $(p=0, q=1)$ of (\ref{204}).
The remaining alternatives 
(modulo $(p,q))$ of this relation
are obtained in a similar way.
\eproof

\section{Main result}
\setcounter{equation}{0}

Let $(X,d,\le)$ be 
a quasi-ordered metric space;
and $T\in \calf(X)$ be a selfmap of $X$.
The basic hypotheses to be considered here are
(a01)-(a03).
\sk

{\bf (A)}
Denote, for each $x,y\in X$
\bit
\item[(c01)]
$M_1(x,y)=d(Tx,Ty)$, $M_2(x,y)=d(x,y)$, $M_3(x,y)=d(x,Tx)$, \\
$M_4(x,y)=d(y,Ty$), $M_5(x,y)=d(x,Ty)$, $M_6(x,y)=d(Tx,y)$, \\
$\calm(x,y)=
(M_1(x,y),M_2(x,y),M_3(x,y),M_4(x,y),M_5(x,y),M_6(x,y))$,\\
$\calm^1(x,y)=(M_2(x,y),M_3(x,y),M_4(x,y),M_5(x,y),M_6(x,y))$.
\eit
Given $F\in \calf(R_+^6,R)$, let us say that $T$ is
{\it $(d,\le;\calm;F)$-contractive}, provided
\bit
\item[(c02)]
$F(\calm(x,y))\le 0$, for all $x,y\in X$ 
with $x\le y$, $x\ne y$.
\eit

The main result of this exposition is

\btheorem \label{t2}
Assume that $T$ is $(d,\le,\calm;F)$-contractive,
for some function $F\in \calf(R_+^6,R)$.
Then, 

{\bf I)}\ \
If $F$ is compatible, almost 2-right-lim-positive
and 4-point-lim-positive,
then $T$ is a Picard operator (modulo $(d,\le)$)

{\bf II)}\ \
If, in addition, $F$ is (3,4)-normal, then
$T$ is a global Picard operator (modulo $(d,\le)$).
\etheorem

\bproof
We first check the final part of the
global Picard property for $T$, by means of  
the extra condition in {\bf II)}.
Let $z_1,z_2\in \Fix(T)$ be such that 
$z_1\le z_2$ and $z_1\ne z_2$. 
By the contractive condition,
$$
F(\rho,\rho,0,0,\rho,\rho)\le 0,\ 
\mbox{where}\
\rho:=d(z_1,z_2)> 0.
$$
This, however, is in contradiction with
$F$ being (3,4)-normal; and our claim follows.
It remains now to establish the 
Picard property, from the conditions in 
{\bf I)}.
Let $x_0\in X(T,\le)$ be arbitrary fixed; 
and put $x_n=T^nx_0$, $n\in N$;
clearly, $(x_n)$ is ascending, by (a03).
Without loss, one may assume that 
$x_n\ne x_{n+1}$, $\forall n$;
note that, in such a case, 
$(r_n:=d(x_n,x_{n+1}); n\ge 0)$
is a sequence in $R_+^0$.

{\bf Step 1.}
Denote for simplicity 
$(s_n:=d(x_{n},x_{n+2}); n\ge 0)$; it is a sequence in 
$R_+$.
By the  contractive condition attached to $(x_{n-1},x_n)$ 
we have
\beq \label{301}
F(r_n,r_{n-1},r_{n-1},r_n,s_{n-1},0)\le 0,\ \ \forall n\ge 1.
\eeq
Combining with the evaluation ($\forall n\ge 1$)
\beq \label{302}
|s_{n-1}-r_{n-1}|=|d(x_{n-1},x_{n+1})-d(x_{n-1},x_{n})|
\le d(x_{n},x_{n+1})=r_n,
\eeq
one gets (via $F$=compatible) that 
\beq \label{303}
\mbox{
($d(x_n,x_{n+1})> 0$,\ $\forall n$,\ and)\  
$d(x_n,x_{n+1})\to 0$ as $n\to \oo$;
}
\eeq
hence, in particular, $(x_n; n\ge 0)$ is $d$-semi-Cauchy.

{\bf Step 2.}
As $F$ is almost 2-right-lim-positive, 
$\Theta:=\Pos(2-right-lim;F)$ is
$(>)$-cofinal in $R_+^0$.
We show that $(x_n; n\ge 0)$ is $d$-Cauchy.
Suppose not; then, by Proposition \ref{p1}, there 
exist a number $b\in \Theta$, a rank $j(b)\ge 0$, 
and a couple
of rank-sequences $(m(j); j\ge 0)$, $(n(j); j\ge 0)$ with
the properties (\ref{201})-(\ref{204}).
By the very definition of $\Theta$, 
$F$ is 2-right-lim-positive at $b$.
On the other hand, (\ref{303}) tells us that
\beq \label{304}
\mbox{
\btab{l}
$(t_3^j:=r_{m(j)}; j\ge 0)$ and $(t_4^j:=r_{n(j)}; j\ge 0)$ \\
are sequences in $R_+^0$ with 
$t_3^j, t_4^j\to 0$ as $j\to \oo$.
\etab
}
\eeq
Moreover, taking (\ref{203}) into account,  yields
\beq \label{305}
\mbox{
\btab{l}
$(t_2^j:=d(x_{m(j)},x_{n(j)}); j\ge 0)$ \\
is a sequence in 
$R_+^0$ with $t_2^j\to b+$ as $j\to \oo$.
\etab
}
\eeq
Finally, by the relation (\ref{204}), one gets
\beq \label{306}
\mbox{
\btab{l}
$(t_1^j:=d(x_{m(j)+1},x_{n(j)+1}); j\ge j(b))$, and \\
$(t_5^j:=d(x_{m(j)},x_{n(j)+1}); j\ge j(b))$, \ $(t_6^j:=d(x_{m(j)+1},x_{n(j)}); j\ge 0)$ \\
are sequences in $R_+^0$ with 
$t_1^j, t_5^j, t_6^j \to b$ as $j\to \oo$.
\etab
}
\eeq
Now, by the first half of (\ref{305}), 
the contractive condition applies to 
$(x_{m(j)},x_{n(j)})$, for all $j\ge 0$; and yields:
$$
F(t_1^j,t_2^j,t_3^j,t_4^j,t_5^j,t_6^j)\le 0,\ \forall j\ge 0.
$$
This gives at once
$\limsup_j F(t_1^j,t_2^j,t_3^j,t_4^j,t_5^j,t_6^j)\le 0$;
hence, $F$ is not 
2-right-lim-positive at $b$; contradiction.

{\bf Step 3.}
As $(x_n; n\ge 0)$ is an ascending $d$-Cauchy sequence, there exists, 
by (a01), some point $x^*\in X$
with $x_n\to x^*$ as $n\to \oo$.
So, if $T$ is $(\le)$-continuous, 
$y_n:=x_{n+1}=Tx_n\to Tx^*$ as $n\to \oo$.
In addition, as $(y_{n}; n\ge 0)$ is a subsequence of 
$(x_n; n\ge 0)$, 
we have $y_n\to x^*$ as $n\to \oo$;
hence (as $d$=metric), $x^*=Tx^*$.
Suppose now that $(\le)$ is self-closed; note that, as a consequence, 
$x_n\le x^*$, $\forall n$.
Two cases may occur.

{\bf Case 3-1.}
There exists a sequence of ranks $(k(i); i\ge 0)$ with $k(i)\to \oo$
as $i\to \oo$ in such a way that 
$x_{k(i)}=x^*$ (hence $x_{k(i)+1}=Tx^*$), for all $i$.
This, and $(x_{k(i)+1}; i\ge 0)$ being a subsequence of 
$(x_n; n\ge 0)$,
gives $x^*\in \Fix(T)$.

{\bf  Case 3-2.}
There exists some rank $h\ge 0$ 
such that 
\bit
\item[(c03)]
$n\ge h$ $\limpl$ $x_n\ne x^*$. 
\eit
Suppose by contradiction that $x^*\ne Tx^*$; i.e.: $b:=d(x^*,Tx^*)> 0$.
From the imposed assumptions, 
$F$ is 4-point-lim-positive at $b$. 
On the other hand, relations (\ref{303})+(c03) 
and our convergence property give (for all $n\ge h$)
\beq \label{307}
\mbox{
\btab{l}
$(t_2^n:=d(x_n,x^*); n\ge h)$, and \\
$(t_3^n:=d(x_n,x_{n+1}); n\ge h)$,
$(t_6^n:=d(x_{n+1},x^*); n\ge h)$ \\
are sequences in $R_+^0$ with
$t_2^n, t_3^n, t_6^n \to 0$ as $n\to \oo$.
\etab
} 
\eeq
Further, the same convergence relation assures us that
$$
\exists n(b)\ge h:\ 0< d(x_n,x^*)< b/2,\ \ \forall n\ge n(b).
$$
Combining with the evaluation
$$
|d(x_n,Tx^*)-b|\le d(x_n,x^*)< b/2,\ \forall n\ge n(b),
$$
we get 
\beq \label{308}
\mbox{
\btab{l}
$(t_1^n:=d(x_{n+1},Tx^*); n\ge n(b))$,
$(t_5^n:=d(x_n,Tx^*); n\ge n(b))$ \\
are sequences in $R_+^0$ with
$t_1^n\to b$, $t_5^n\to b$ if $n\to \oo$.
\etab
}
\eeq
The contractive condition
applies to $(x_n,x^*)$ (for $n\ge 0$); and yields 
$$
F(t_1^n,t_2^n,t_3^n,t_4^n,t_5^n,t_6^n)\le 0,\ 
\forall n\ge n(b);
$$
where $(t_4^n=b; n\ge 0)$.
This gives at once 
$\limsup_n F(t_1^n,t_2^n,t_3^n,t_4^n,t_5^n,t_6^n)\le 0$;
hence, $F$ is not 4-point-lim-positive at $b$;
contradiction. 
So, necessarily, $x^*\in \Fix(T)$; and the proof is complete.
\eproof

\section{Explicit versions}
\setcounter{equation}{0}

In the following, we show that 
a certain "explicit" version of this result
yields a quasi-order extension of the 2008 one in
Agarwal et al \cite{agarwal-el-gebeily-o-regan-2008}.
\sk

Given $k\ge 1$, let us say that $G\in \calf(R_+^k,R_+)$ 
is a {\it $k$-altering function} provided
\bit
\item[(d01)]
$G$ is continuous and increasing in all variables
\item[(d02)]
$G$ is reflexive-sufficient:\ 
$G(t_1,...,t_k)=0$ iff $t_1=...=t_k=0$.
\eit
The class of all these will be denoted as 
$\calf(alt)(R_+^k,R_+)$.
When $k=4$, a basic example of 
4-altering function may be constructed as
\bit
\item[(d03)]
$L(t_1,t_2,t_3,t_4)=\max\{t_1,t_2,t_3,t_4\}$,
$(t_1,t_2,t_3,t_4)\in R_+^4$.
\eit
Let $L^*\in \calf(R_+^5,R_+)$ be the associated function
\bit
\item[(d04)]
$L^*(t_1,t_2,t_3,t_4,t_5)=L(t_1,t_2,t_3,(1/2)(t_4+t_5))$,
$(t_1,t_2,t_3,t_4,t_5)\in R_+^5$;
\eit
note that $L^*$ is a 5-altering function,
as it can be directly seen. 

Let $(X,d,\le)$ be 
a quasi-ordered metric space;
and $T\in \calf(X)$ be a selfmap of $X$.
The basic hypotheses to be considered here are
again (a01)-(a03).
For each $x,y\in X$,
let $(M_i(x,y); 1\le i\le 6)$,
$\calm(x,y)$ and $\calm^1(x,y)$ 
be the ones of (c01).
Given $\psi\in \calf(R_+)$, 
let us say that $T$ is
{\it $(d,\le;\calm;\psi)$-contractive} if 
\bit
\item[(d05)]
$d(Tx,Ty)\le \psi(L^*(\calm^1(x,y)))$, 
$\forall x,y\in X$, $x\le y$, $x\ne y$.
\eit
Note that the introduced convention amounts to saying that 
$T$ is $(d,\le;\calm;F)$-contractive, where
$F\in\calf(R_+^6,R)$ is introduced as
\bit
\item[(d06)]
$F(t_1,...,t_6)=t_1-\psi(L^*(t_2,...,t_6)),\ 
(t_1,...,t_6)\in R_+^6$.
\eit
We want to determine under which conditions about $\psi$ 
is Theorem \ref{t2} applicable to $(X,\le,d)$ and the function $F$.
\sk

{\bf (A)}
Given $\psi\in \calf(re)(R_+)$,  call it  {\it compatible}, when 
\bit
\item[(d07)]
for each sequence $(r_n; n\ge 0)$ in $R_+^0$, with
$r_n\le \psi(r_{n-1})$, $\forall n\ge 1$,
we must have $r_n\to 0$.
\eit

To get sufficient conditions for
such a property,  denote, for each $s> 0$
\bit
\item[(d08)]
$P(s):=\limsup_{t\to s+} \psi(t)$,\ \  $Q(s)=\max\{\psi(s),P(s)\}$.
\eit
Clearly, by the regressive property of $\psi$, we must have
$P(s)\le s$; wherefrom, 
$\psi(s)\le Q(s)\le s$ (for all $s> 0$).
Note that neither of the inequalities above is
strict (at some $s> 0$), in general. 
For, if $\psi$ is 
continuous from the right at $s$, we have $Q(s)=\psi(s)$;
so, the former of these inequalities is non-strict.
On the other hand, if $\psi$ is increasing 
(at least on an open interval containing $s$), 
$Q(s)=\max\{\psi(s),\psi(s+0)\}$; hence, 
whenever $\psi(s+0)=s$, the latter of these inequalities is non-strict.

Call $\psi\in \calf(re)(R_+)$, {\it Boyd-Wong-admissible} 
at $s> 0$, when
$P(s)< s$ (or, equivalently: $Q(s)< s$).
If this holds for all $s> 0$, then we shall say that 
$\psi$ is {\it Boyd-Wong admissible}.

\blemma \label{le1}
For each $\psi\in \calf(re)(R_+)$: Boyd-Wong admissible
$\limpl$ compatible.
\elemma

\bproof
Let $(r_n; n\ge 0)$ be a sequence in $R_+^0$ with 
$r_n\le \psi(r_{n-1})$, $\forall n\ge 1$.
As $\psi\in \calf(re)(R_+)$, $(r_n)$ is strictly
descending in $R_+$; hence, $r:=\lim_n r_n$ exists
in $R_+$ and [$r_n> r$, $\forall n$].
We have (again via $\psi\in \calf(re)(R_+)$) 
$r_n\le \psi(r_{n-1})< r_{n-1}$, $\forall n\ge 1$.
This, along with $r_n \to r$ as $n\to \oo$, yields 
$\lim_n \psi(r_n)=r$; wherefrom $P(r)=r$;
contradiction. Hence, $r=0$, as desired.
\eproof

Having these precise, call $\psi\in \calf(re)(R_+)$,
{\it almost Boyd-Wong admissible} when 
for each $\veps> 0$ there exists $s\in ]0,\veps[$ 
with $Q(s)< s$.

\bprop \label{p2}
Let the function $\psi\in \calf(re)(R_+)$ be compatible
and almost Boyd-Wong admissible.
Then, the function $F$ given by (d06) is 
compatible, almost 2-right-positive,
4-point-lim positive, and (3,4)-normal.
\eprop

\bproof
The argument will be divided into several steps.

{\bf Part 1}\ ($F$ is compatible).
Let $(r_n)\subset R_+^0$, $(s_n)\subset R_+$ be sequences fulfilling
\bit
\item[(d09)]
$F(r_n,r_{n-1},r_{n-1},r_n,s_{n-1},0)\le 0$
and $|s_{n-1}-r_{n-1}|\le r_{n}$, $\forall n\ge 1$.
\eit
From $s_{n-1}\le r_{n-1}+r_n\le 2\max\{r_{n-1},r_n\}$, $\forall n\ge 1$, 
we have
$$
L^*(r_{n-1},r_{n-1},r_n,s_{n-1},0)= \max\{r_{n-1},r_n\}, \forall n\ge 1; 
$$
so that, the above inequality becomes 
$r_n\le \psi(\max\{r_{n-1},r_n\})$, $\forall n\ge 1$.
This, via $\psi\in \calf(re)(R_+)$, gives $r_n\le\psi(r_{n-1})$,
$\forall n\ge 1$; wherefrom (as $\psi$ is compatible) 
$r_n\to 0$ as $n\to \oo$; and the claim follows.

{\bf Part 2}\ ($F$ is (3,4)-normal).
Let $r> 0$ be arbitrary fixed. By definition,
$$
F(r,r,0,0,r,r)=r-\psi(L^*(r,0,0,r,r))=
r-\psi(r)> 0;
$$
and, from this, we are done.

{\bf Part 3}\ ($F$ is almost 2-right-lim-positive).
As $\psi$ is almost Boyd-Wong admissible,
for each $\veps> 0$ there exists $r\in ]0,\veps[$
with $Q(r)< r$. 
We show that the function $F$ defined above is 
2-right-lim-positive at $r$.
Let $(t_i^n; n\ge 0)$, $i\in \{1,2,3,4,5,6\}$, be  
sequences in $R_+^0$ with (as $n\to \oo$) 
\bit
\item[]
$t_i^n\to r$, $i\in \{1,5,6\}$;\ 
$t_i^n\to 0$, $i\in \{3,4\}$;\
$t_2^n\to r+$.
\eit
By definition, this yields
$\ga_n:=L^*(t_2^n,t_3^n,t_4^n,t_5^n,t_6^n)\ge t_2^n> r$, $\forall n$;
and, as $L^*$ is continuous in its variables,
$\ga_n\to L^*(r,0,0,r,r)=r$ as $n\to \oo$;
hence, summing up, $\ga_n\to r+$ as $n\to \oo$.
As a consequence,
$$
\limsup_n F(t_1^n,t_2^n,t_3^n,t_4^n,t_5^n,t_6^n)\ge
r-\liminf_n\psi(\ga_n)\ge r-Q(r)> 0;
$$
and this proves our assertion.

{\bf Part 4}\ ($F$ is 4-point--lim-positive).
Let $r> 0$ be arbitrary fixed.
We have to show that $F$ is 4-point-lim-positive 
at $r$.
Let $(t_i^n)$, $i\in \{1,2,3,4,5,6\}$, be  
sequences in $R_+^0$ 
with $t_4^n=r$, $\forall n$; and
(as $n\to \oo$)
\bit
\item[] 
$t_i^n\to r$, $i\in \{1,5\}$;\ 
$t_i^n\to 0$, $i\in \{2,3,6\}$. 
\eit
There exists some 
rank $n(r)$ in such a way that ($\forall n\ge n(r)$)
$t_i^n< 3r/2$, $i\in \{1,5\}$, and $t_i^n< r/2$, $i\in \{2,3,6\}$. 
Combining with the choice of $(t_4^n)$, yields
$$
\ga_n:=L^*(t_2^n,t_3^n,t_4^n,t_5^n,t_6^n)=r,\ 
\forall n\ge n(r);
$$
wherefrom 
$\limsup_n F(t_1^n,t_2^n,t_3^n,t_4^n,t_5^n,t_6^n)=r-\psi(r)> 0$.
\eproof

Now, by simply combining the obtained fact  
with Theorem \ref{t2}, one gets
the following practical statement.
(As before, the basic hypotheses (a01)-(a03)
prevail).

\btheorem \label{t3}
Suppose that $T$ is $(d,\le;\calm;\psi)$-contractive,
for some compatible almost Boyd-Wong admissible
$\psi\in \calf(re)(R_+)$.
Then, $T$ is a Picard operator (modulo $(d,\le)$).
\etheorem
\sk

{\bf (B)}
Let us now give some particular cases of
this result, with a practical finality.

{\bf B-1)}
Suppose that $\psi\in \calf(in,re)(R_+)$ is 
{\it Matkowski-admissible} (cf. Section 1).
Let $(r_n)$ be a sequence in $R_+^0$ with 
$r_n\le \psi(r_{n-1})$, $\forall n\ge 1$.
As $\psi$ is increasing, this yields 
$r_n\le \psi^n(r_0)$, $\forall n$; wherefrom 
$r_n\to 0$; hence $\psi$  is  compatible.
On the other hand, let $\Ga:=\Ga_\psi$ stand for
the (at most denumerable) subset of all $r> 0$
where $\psi$ is discontinuous.
Each $r> 0$ not belonging to $\Ga$
is a (bilateral) continuity point of $\psi$;
and then (as $\psi(r)=Q(r)$), $\psi$ is 
Boyd-Wong admissible at $r$.
Summing up, $\psi$ is compatible and
almost Boyd-Wong admissible.
The corresponding version of Theorem \ref{t3} 
under this choice of $\psi$ is just 
the 2008 fixed point statement in 
Agarwal et al \cite{agarwal-el-gebeily-o-regan-2008};
see also 
O'Regan and Petru\c{s}el \cite{o-regan-petrusel-2008}.

{\bf B-2)}
Suppose now that
$\psi\in \calf(re)(R_+)$ is Boyd-Wong-admissible.
Clearly, $\psi$ is compatible, by Lemma \ref{le1};
moreover (by definition), $\psi$ is 
almost Boyd-Wong admissible.
The corresponding version of Theorem \ref{t3}
under this choice of $\psi$ is a counterpart
of the above cited 2008 result in 
Agarwal et al \cite{agarwal-el-gebeily-o-regan-2008};
but it cannot be reduced to it.
In particular, when $(\le)=X\times X$,
the same variant includes 
the fixed point result in 
Boyd and Wong \cite{boyd-wong-1969};
as well as (when $\psi$ is linear), the result in
Hardy and Rogers \cite{hardy-rogers-1973}.

\section{Global aspects}
\setcounter{equation}{0}

In the following, a certain "global" version 
of the main result is given. 
As before, $(X,\le,d)$ is a quasi-ordered  metric space; 
and $T\in \calf(X)$ is a selfmap of $X$.
\sk

Let the function $F\in\calf(R_+^6,R)$ be compatible
[in the sense of (b01)]. 
For an application of Theorem \ref{t2} it will suffice
that $F$ be (in addition)
almost 2-right-lim-positive, 4-point-lim-positive, and
(eventually) normal.
We shall try to assure this under the global condition
\bit
\item[(e01)]
$F$ is lower semicontinuous (in short: lsc) on $R_+^6$:  \\
$\liminf_n F(t_1^n,...,t_6^n)\ge F(a_1,...,a_6)$,
whenever $t_i^n\to a_i$, $i\in \{1,...,6\}$.
\eit
Note that, in such a case, the lim-positive conditions 
are obtainable from 
\bit
\item[(e02)]
$F(r,r,0,0,r,r)> 0$, $F(r,0,0,r,r,0)> 0$, $\forall r> 0$;
\eit
referred to as: $F$ is {\it (3,4)-normal}
and {\it (2,3,6)-normal}, respectively;
the former of these is just condition (b02).

An application of Theorem \ref{t2} yields
the following practical result.
(The basic hypotheses 
to be considered here are
again (a01)-(a03)).

\btheorem \label{t4}
Assume that 
$T$ is $(d,\le,F)$-contractive,
for some compatible lsc $F\in \calf(R_+^6,R)$
which is both (3,4)-normal and (2,3,6)-normal.
Then, 
$T$ is a global Picard operator (modulo $(d,\le)$).
\etheorem

The following particular case is of interest.
Assume that (in addition to (e01))
the global condition holds
\bit
\item[(e03)]
$F$ is $(2,...,6)$-decreasing:
$F(t_1,.)$ is decreasing, $\forall t_1\in R_+$.
\eit
Then, the compatibility condition (b01) 
is deductible from:
\bit
\item[(e04)]
$F$ is almost-compatible: for each sequence $(r_n)$, with \\
$F(r_n,r_{n-1},r_{n-1},r_n,r_n+r_{n-1},0)\le 0$, $\forall n\ge 1$,
we must have $r_n\to 0$.
\eit
In particular this last condition is deductible (via Lemma \ref{le1}) from
\bit
\item[(e05)]
$F$ is $\psi$-compatible
($F(u,v,v,u,u+v,0)\le 0$ $\limpl$ $u\le \psi(v)$) \\
for some admissible function $\psi\in \calf(re)(R_+)$.
\eit
This is just the main result in
Altun and Simsek \cite{altun-simsek-2010};
obtained (under a different approach)
with (e01) being substituted by a
continuity assumption about $F$.
[In fact, the authors' argument cannot be
entirely acceptable; for, e.g., the implication 
(3.8) $\limpl$ (3.9) in that paper; i.e.
$$
[r_n\le \psi(r_{n-1}), \forall n\ge 1] \limpl 
r_n\le \psi^n(r_0), \forall n
$$
is not true unless $\psi$ is increasing.
The same remark is valid for Lemma 3.3 in that paper;
we do not give details].

Now, technically speaking, condition (e05)
was introduced so as to be applicable to functions
$F$ like in (d06), where 
the admissible $\psi\in \calf(re)(R_+)$ is 
either increasing or continuous.
In the former case, $F$ is
$(2,...,6)$-decreasing; but, not in general lsc.
In the latter case, $F$ is neither 
lsc nor $(2,...,6)$-decreasing.
As a consequence of this,
neither Theorem \ref{t3} nor Theorem \ref{t4}
are deductible from the above result.
[Note that the second half of this conclusion 
is in contradiction to the authors' claim,
expressed via Example 2.3 of their paper].
Further aspects 
may be found in 
Popa and Mocanu \cite{popa-mocanu-2009};
see also 
Vetro and Vetro \cite{vetro-vetro-2013}.

\section{Old approach (1976)}
\setcounter{equation}{0}

In the following, a summary of the 1976 results in 
Turinici \cite{turinici-1976} 
is being sketched, for historical reasons
(explained at the beginning).
\sk 

Let $S \ne \es$ be a nonempty set and 
$P$, some nonempty proper subset of it ($\es \ne P \subset S$) .
Denote by $(S6)$ the class of all functions 
$F:R_+^6 \to S$; 
and by $(P6)$, the subclass of all $F\in (S6)$ 
satisfying the global conditions
\bit
\item[(f01)] 
$w> 0$ $\impl$ $F(w,w,0,0,w,w)\in S\sm P$
\item[(f02)] 
$u,v> 0$, $p\le u+v$, $F(u,v,v,u,p,0)\in P$ $\impl$ $u\le v$     
\eit
as well as the local conditions: $\forall r> 0$, $\exists a(r)\in ]0,r[$ such that 
\bit
\item[(f03)]
$u,v\in [r,r+a(r)[$, $u\le v$, $p\le u+v$ $\impl$ 
$F(u,v,v,u,p,0) \in S\sm P$
\item[(f04)]
$t,p,q\in]r-a(r),r+a(r)[,$ $u\in [r,r+a(r)[$, $v,w\in ]0,a(r)[$ $\impl$ \\
$F(t,u,v,w,p,q)\in S\sm P$
\item[(f05)] 
$t,p \in ]r-a(r),r+a(r)[$, $u,v,q \in ]0,a(r)[$ $\impl$  
$F(t,u,v,r,p,q) \in S\sm P$.
\eit

Having these precise,
let $(X,d)$ be a complete metric space;
and $T:X\to X$, a selfmap of $X$.
Given $F\in (S6)$, 
we say that $T$ is a  {\it 6-implicit contraction mapping}
(abbreviated: 6-icm) with respect to it, provided
\bit
\item[(f06)]
$F(d(Tx, Ty),d(x,y),d(x,Tx),d(y,Ty),d(x,Ty),d(y,Tx))\in P$ \\
for all $x,y\in  X$ with $Tx \ne Ty$.
\eit

Our main result can be stated as follows.

\btheorem \label{t5}
Suppose that $T:X \to X$ is a 6-icm with respect to some $F\in (P6)$.
Then the following conclusions hold:
\beq  \label{601}
\mbox{
$T$ has a unique fixed point $z \in X$
}
\eeq
\beq \label{602}
\mbox{
$T^nx \to z$,\ as  $n \to \oo$,\ $\forall x\in X$.
}
\eeq
\etheorem 

\bproof
First, we prove the uniqueness of the fixed point of $T$. 
Let $z_1,z_2\in X$ be such that $z_1 = Tz_1$, $z_2 = Tz_2$, 
$z_1 \ne z_2$. 
From (f06) and (f01), we obtain
$$
F(d(z_1,z_2),d(z_1,z_2),0,0,d(z_1,z_2),d(z_1,z_2)) 
\in P\cap (S\sm P)=\es,
$$
contradiction; therefore, $z_1=z_2$. 
Now we prove the existence. 
Take any $x_0 \in X$ and consider the sequence $\{x_n:=T^nx_0; n \ge 0\}$. 
If $x_n=x_{n+1}$ for some $n$, the conclusion follows. 
Assume that $x_n \ne x_{n+1}$, $\forall n\ge 0$. 
From (f06), 
\beq \label{603}
\barr{l}
F(d(x_n,x_{n+1}), d(x_{n-1},x_n), d(x_{n-1},x_n), \\
d(x_n,x_{n+1}), d(x_{n-1},x_{n+1}),0)\in P,\  \forall n\ge 1.
\earr
\eeq
On the other hand, the triangle inequality gives ($\forall n\ge 1$)
\beq \label{604}
d(x_{n-1},x_{n+1})\le d(x_{n-1},x_n) + d(x_n,x_{n+1}).
\eeq
From (\ref{603}), (\ref{604}) and (f02) we obtain ($\forall n\ge 1$):
\beq \label{605}
d(x_n,x_{n+1}) \le d(x_{n-1}, x_n);
\eeq
i.e., the sequence $\{r_n:=d(x_nx,x_{n+1}); n\ge 0\}$ decreases. 
Let $r=\lim_n d(x_n, x_{n+1})$ and assume that $r> 0$. 
One can find some rank $n(r)\ge 1$ such that
\beq \label{606}
n\ge n(r) \impl d(x_{n-1},x_n) \in [r,r+a(r)[.
\eeq
Taking into account (\ref{604})--(\ref{606}) and (f03), we have
for all $n\ge n(r)$
\beq \label{607}
\barr{l}
F(d(x_n,x_{n+1}),d(x_{n-1},x_n),d(x_{n-1},x_n), \\
d(x_n,x_{n+1}),
d(x_{n-1},x_{n+1}),0) \in S\sm P;
\earr
\eeq
which contradicts (\ref{603}) for $n\ge n(r)$. Therefore $r=0$.
Suppose that $\{x_n; n \ge 0\}$ is not a Cauchy sequence. 
Then there exist $\veps > 0$ and two sequences of natural numbers 
$\{m(j); j\ge 0\}$ and $\{n(j); j\ge 0\}$, $m(j)< n(j)$, $m(j)\to \oo$ as $j\to \oo$ 
such that 
$d(x_{m(j)},x_{n(j)})\ge \veps$, while 
$d(x_{m(j)}, x_{n(j)-1})< \veps$, $\forall j\ge 0$. 
For the sake of simplicity, we shall write $m, n$, instead of $m(j), n(j)$, respectively.
As $d(x_k,x_{k+1}) \to 0$ as $k \to \oo$, we can find $j(\veps) \in N$ such that
($\forall j\ge j(\veps$))
\beq \label{608}
0< d(x_{n-1}, x_n)\le d(x_m, x_{m+1})< (1/3)a(\veps)< a(\veps)< \veps.
\eeq
On the other hand, from the triangle inequality we have ($\forall j \ge 0$)
\beq \label{609}
\barr{l}
d(x_m,x_n)-d(x_m,x_{m+1})-d(x_n,x_{n+1})\le \\ d(x_{m+1},x_{n+1}) \le 
d(x_m,x_n)+d(x_m,x_{m+1})+d(x_n,x_{n+1}),
\earr
\eeq
\beq \label{610}
d(x_m, x_n)\le d(x_m, x_{n-1})+d(x_{n-1}, x_n),
\eeq
\beq \label{611} 
\barr{l}
d(x_m,x_n)-d(x_n,x_{n+1}) \le d(x_m, x_{n+1})\le  d(x_m,x_n)+d(x_n,x_{n+1}),
\earr
\eeq
\beq \label{612}
\barr{l}
d(x_m,x_n)-d(x_m,x_{m+1})\le d(x_n,x_{m+1}) \le  
d(x_m, x_n)+d(x_m,x_{m+1}). 
\earr
\eeq
From (\ref{608})--(\ref{612}) it easily follows, $\forall j \ge j(\veps)$
\beq \label{613}
\barr{l}
d(x_{m+1},x_{n+1}),d(x_m,x_{n+1}),d(x_n,x_{m+1})\in \\
]\veps-a(\veps), \veps+a(\veps)[, \  \
d(x_m,x_n)\in [\veps, \veps+a(\veps)[.
\earr
\eeq
Now, (f06), (\ref{608}), (\ref{613}) and (f04) give us (for all  $j \ge j(\veps)$)
$$ 
\barr{l}
F(d(x_{m+1}, x_{n+1}),d(x_m, x_n),d(x_m,x_{m+1}), \\
d(x_n,x_{n+1}), d(x_m,x_{n+1}), d(x_n,x_{m+1})) \in 
P\cap (S\sm P)=\es,
\earr
$$
a contradiction. Therefore, $\{x_n; n \ge 0\}$ is a Cauchy sequence. 
Since $(X, d)$ is complete, $x_n \to z$, for some $z\in  X$. 
We have two possibilities:

{\bf i)}\
There exists a sequence of natural numbers $\{k(n); n\ge 0\}$, 
$k(n)\to \oo$ as $n \to \oo$,
such that $x_{k(n)}= z$. Then, $x_{k(n)+1}=Tz$. 
Letting $n$ tends to infinity and using the fact that 
$\{x_{k(n)+1}; n \ge 0\}$ is a subsequence of $\{x_n; n \ge 0\}$ we get $z= Tz$.

{\bf ii)}\ 
There exists $n_0 \in N$ such that 
$n\ge n_0$ $\impl$ $x_n\ne  z$. 
Suppose that $z \ne Tz$; then $r=d(z, Tz)> 0$. 
We can find $n(r)\in N$, such that, $\forall n\ge n(r)$
\beq \label{614} 
0< d(x_n,x_{n+1}), d(x_n,z)< (1/3)a(r)< a(r)< r.
\eeq
On the other hand, from the triangle inequality, we have 
$$
r-d(x_n,z)\le d(x_n,Tz) \le r+d(x_n, z),\ \forall n\ge 0;
$$
so that (from (\ref{614}))
\beq \label{615} 
d(x_n,Tz) \in ]r-a(r),r+a(r)[, \forall n\ge n(r).
\eeq
Now, (f06), (\ref{614}), (\ref{615}) and (f05) give us for $n\ge n(r)$
$$
\barr{l}
F(d(x_{n+1},Tz),d(x_n,z),d(x_n,x_{n+1}), \\
r,d(x_n,Tz),d(z, x_{n+1}))\in P\cap (S\sm P)=\es,
\earr
$$
a contradiction. Therefore, $z=Tz$, which completes the proof.
\eproof

\brem \label{r1}
\rm

In the original paper, 
the extra requirement below is being added
\bit
\item[(f07)]
$t,u> 0$, $v,w,p,q\ge 0$ $\impl$ \\
$F(t,u,v,w,p,q)=F(t,u,w,v,p,q)=F(t,u,v,w,q,p)$.
\eit
But, evidently, this condition 
(imposed for symmetry reasons) is superfluous.
\erem

\brem \label{r2}
\rm

In particular, letting $S=R$, $P=R_+$,
conditions of Theorem \ref{t5} are comparable with the 
standard ones.
\erem

Now, Theorem \ref{t5} is a partial extension of a result due to 
Hardy and Rogers  \cite{hardy-rogers-1973}.
On the other hand, if 
$(t,u,v,w,p,q)\mapsto F(t,u,v,w,p,q)$ does not depend on its 
last two variables, 
the corresponding form of this result 
extends the ones in  
Reich  \cite{reich-1972}
and 
Turinici \cite{turinici-1974}.
Finally, when
$(t,u,v,w,p,q)\mapsto F(t,u,v,w,p,q)$ does not depend on its 
last four variables, 
Theorem \ref{t5} reduces to the 
fixed point statement in 
Boyd and Wong \cite{boyd-wong-1969}.
Note that this extension
assured by Theorem \ref{t5} is rather 
different from the one in Theorem \ref{t2}.
So, it would be natural asking whether 
a common version of  both these results is possible.
Further aspects will be delineated elsewhere.


\end{document}